\documentclass{elsart}
\usepackage{amssymb,amsfonts}
\newcommand{\trace}{\mathop{\rm Tr}\nolimits}

\newcommand{\twomat}[4]{\left(\begin{array}{cc}#1&#2\\#3&#4\end{array}\right)}
\newcommand{\twovec}[2]{\left(\begin{array}{c}#1\\#2\end{array}\right)}
\newcommand{\twovect}[2]{\left(#1\quad#2\right)}

\newcommand{\cA}{{\mathcal A}}

\newcommand{\cH}{{\mathcal H}}

\newcommand{\bR}{{\mathbb{R}}}

\newcommand{\id}{\mathbb{I}}

\newcommand{\be}{\begin{equation}}
\newcommand{\ee}{\end{equation}}
\newcommand{\bea}{\begin{eqnarray}}
\newcommand{\eea}{\end{eqnarray}}
\newcommand{\beas}{\begin{eqnarray*}}
\newcommand{\eeas}{\end{eqnarray*}}

\newtheorem{theorem}{Theorem}

\newtheorem{definition}{Definition}

\newtheorem{conjecture}{Conjecture}
\newcount\minute
\newcount\hour
\def\currenttime{%
    \minute\time
    \hour\minute
    \divide\hour60
    \the\hour:\multiply\hour60\advance\minute-\hour\the\minute}
\begin{document}
\begin{frontmatter}
\title{On a Block Matrix Inequality quantifying the Monogamy of the Negativity of Entanglement}
\author{Koenraad M.R.\ Audenaert}
\address{
Department of Mathematics,
Royal Holloway University of London, \\
Egham TW20 0EX, United Kingdom \\[1mm]
Department of Physics and Astronomy, University of Ghent, \\
S9, Krijgslaan 281, B-9000 Ghent, Belgium}
\date{\today, \currenttime}
%
\begin{abstract}
We convert a conjectured inequality from quantum information theory, due to He and Vidal, into a block matrix inequality and prove a special case.
Given $n$ matrices $A_i$, $i=1,\ldots,n$, of the same size, let $Z_1$ and $Z_2$ be the block matrices $Z_1:=(A_jA_i^*)_{i,j=1}^n$ and
$Z_2:=(A_j^*A_i)_{i,j=1}^n$. Then the conjectured inequality is
\[
\left(||Z_1||_1-\trace Z_1\right)^2 + \left(||Z_2||_1-\trace Z_2\right)^2
\le \left(\sum_{i\neq j} ||A_i||_2 ||A_j||_2\right)^2.
\]
We prove this inequality for the already challenging case $n=2$ with $A_1=\id$.
\end{abstract}
\end{frontmatter}
\section{Introduction}
Quantum Information Theory (QIT), a recent physical theory combining concepts of information theory with quantum mechanics,
has proven to be a rich source of challenging matrix analysis problems \cite{cargese,AK}.
In this paper one such problem is presented and some progress towards its resolution is reported.

The problem is as follows.
Consider a set of $n$ given general $n_1\times n_2$ matrices $A_i$, and with them form the two $n\times n$ block matrices
\[
Z_1:=(A_j A_i^*)_{i,j=1}^n=
\left(
\begin{array}{ccc}
A_1A_1^* & A_2A_1^* & \ldots \\
A_1A_2^* & A_2A_2^* & \ldots \\
\vdots & \vdots &
\end{array}
\right)
\]
and
\[
Z_2:=(A_j^* A_i)_{i,j=1}^n=
\left(
\begin{array}{ccc}
A_1^*A_1 & A_2^*A_1 & \ldots \\
A_1^*A_2 & A_2^*A_2 & \ldots \\
\vdots & \vdots &
\end{array}
\right).
\]
These two matrices are Hermitian, but not in general positive semidefinite.
Thus, the quantities $||Z_i||_1-\trace Z_i$ are not necessarily zero.
We wish to investigate whether the following inequality holds:
\[
\left(||Z_1||_1-\trace Z_1\right)^2 + \left(||Z_2||_1-\trace Z_2\right)^2
\le \left(\sum_{i\neq j} ||A_i||_2 ||A_j||_2\right)^2.
\]
Here $||\cdot||_1$ and $||\cdot||_2$ denote the trace norm (Schatten 1-norm) and Frobenius norm (Schatten 2-norm), respectively.

This inequality is the block matrix formulation of an equivalent inequality in QIT, conjectured recently by He and Vidal in \cite{HV},
regarding the so-called `monogamy of the negativity of entanglement'.
For the benefit of readers who are not familiar with the QIT jargon we give a brief presentation of this conjecture in Section \ref{sec:not},
where we also explain the QIT concepts of partial trace, partial transpose and negativity, in terms of which the conjecture is expressed.
In Section \ref{sec:convert} we show how the He-Vidal conjecture can be expressed in terms of block matrices, yielding the 
abovementioned inequality, which is (\ref{ineq4}) in that Section. Both Section \ref{sec:not} and \ref{sec:convert} can be skipped by readers who
are not interested in the QIT-background of the problem.

In our opinion, proving this inequality is a very hard problem, and we have only succeeded in proving a very special case.
Namely, we have only been able to prove the case $n=2$, where there are only two matrices $A_1$ and $A_2$, and where in addition we also
require $A_1$ to be the identity matrix. This proof is presented in Section \ref{sec:proof}.

We end this introduction by recalling some of the notations we will use.
The modulus of a matrix $X$ will be denoted as $|X|$, and is given by $(X^*X)^{1/2}$. Any Hermitian matrix can be decomposed as a difference
of its positive and negative part: $X=X_+ - X_-$, with $X_\pm:=(|X|\pm X)/2$. This is the so-called Jordan decomposition.
The Schatten $q$-norm of a matrix, for $q\ge 1$ is denoted as $||X||_q$ and is defined as $||X||_q := (\trace |X|^q)^{1/q}$. 
The trace norm is just the Schatten 1-norm, $||X||_1 = \trace X$, and the Frobenius norm is the Schatten 2-norm.
We will also need the quantity $||X||_q$ for $0<q<1$, which is no longer a norm but a quasi-norm.
Finally, we denote the eigenvalues of a Hermitian matrix, sorted either in non-increasing or non-decreasing order as $\lambda^\downarrow_j$
and $\lambda^\uparrow_j$, respectively.
\section{The He-Vidal Conjecture in Quantum Information Theory \label{sec:not}}
Let us begin with highlighting some of the main mathematical features about QIT in general, and the problem in particular.
We will be very brief and refer to \cite{cargese} for a more in-depth discussion.

For convenience (and for interfacing with the QIT part of the readership)
we will use Dirac notation for vectors up until Section \ref{sec:convert}, after which it is no longer needed.
A general vector of a Hilbert space $\cH$ will be denoted as $|\psi\rangle$, where the symbol $\psi$ is merely a label;
the greek letters $\psi$ and $\phi$ are typically reserved for this purpose.
The Hermitian conjugate of the vector $|\psi\rangle$ is denoted by $\langle \psi|$.
The inner product between two such vectors is $\langle\psi|\phi\rangle$, whereas the outer product
is $|\psi\rangle\langle\phi|$ (for $\psi\otimes\phi^*$).
The tensor product of two vectors $|\phi\rangle$ and $|\psi\rangle$ is denoted $|\phi\rangle\otimes|\psi\rangle$
or $|\phi\rangle|\psi\rangle$ for short.

The elements of an orthonormal basis for a finite-dimensional Hilbert space $\cH$ (with dimension $d$)
will be denoted by $|i\rangle$, with $i=1,\ldots,d$.
To distinguish between basis vectors of bases for several Hilbert spaces, various roman letters will be used ($i$, $j$, $k$)
serving the dual purpose of label and of index.

One of the main tenets of quantum mechanics is that the state of any quantum-mechanical system (a set of atoms, say)
is completely described by a complex, normalised vector in some Hilbert space $\cH$. This vector
is the  wavevector or wavefunction known from any introductory quantum mechanics course.
In this paper we will only
study finite-dimensional quantum systems, having Hilbert spaces of finite dimension $d$.
However, in the experimentally more relevant case that only part of the quantum system is accessible to experiment
(see \cite{cargese} for a detailed explanation of what that means)
one has to resort to another description whereby a state is represented by a positive semidefinite
$d\times d$ matrix with trace 1, known as a density matrix and usually denoted by a greek letter, such as $\rho$.
The two descriptions coincide when the whole quantum system is experimentally accessible.
In that case the density matrix has rank 1 and can be written as $\rho=|\psi\rangle\langle\psi|$,
where $|\psi\rangle$ is exactly the wavevector of the quantum system.

The problem that we wish to study here involves quantum systems composed of three subsystems (`parties')
labeled $A$, $B$ and $C$, each of which has an associated Hilbert space $\cH_A$, $\cH_B$ and $\cH_C$
of dimension $d_A$, $d_B$ and $d_C$, respectively. The Hilbert space
$\cH$ of the entire system is the tensor product $\cH=\cH_A\otimes\cH_B\otimes \cH_C$ and has dimension $d=d_Ad_Bd_C$.
Let us denote the orthonormal bases of these Hilbert spaces by $\{|i\rangle\}_{i=1}^{d_A}$ (for $\cH_A$),
$\{|j\rangle\}_{j=1}^{d_B}$ (for $\cH_B$) and $\{|k\rangle\}_{k=1}^{d_C}$ (for $\cH_C$).
The most natural basis of $\cH$ (at least from the quantum information perspective) is the tensor product
of these three bases: $\{|i\rangle|j\rangle|k\rangle\}_{i,j,k}$, or even shorter $\{|ijk\rangle\}_{i,j,k}$.
Given a density matrix $\rho$ on $\cH$, its matrix elements in the supplied basis are
$\langle ijk|\rho|i'j'k'\rangle$.

Two important mathematical operations that are prominent in QIT are relevant for this paper:
the partial trace and the partial transpose. Crudely speaking, these are the familiar operations of trace and transpose but
acting on one subsystem only.

We will need the partial traces $\trace_B$ and $\trace_C$.
Given a density matrix $\rho$ on $\cH=\cH_A\otimes\cH_B\otimes \cH_C$,
these operations yield density matrices on $\cH_A\otimes\cH_C$
and $\cH_A\otimes\cH_B$, respectively.
The matrix elements of $\trace_B\rho$ are given by
\[
\langle ik|\trace_B \rho|i'k'\rangle = \sum_{j=1}^{d_B} \langle ijk|\rho|i'jk'\rangle
\]
and those of $\trace_C\rho$ by
\[
\langle ij|\trace_C \rho|i'j'\rangle = \sum_{k=1}^{d_C} \langle ijk|\rho|i'j'k\rangle.
\]

The partial transpose (with respect to party $A$) of a density matrix $\rho$ on $\cH$, denoted $\rho^\Gamma$, is
a Hermitian matrix on $\cH$ (not necessarily positive semidefinite anymore) with matrix elements
\[
\langle ijk|\rho^\Gamma|i'j'k'\rangle = \langle i'jk|\rho|ij'k'\rangle.
\]
If $\rho$ is a tensor product of two density matrices, $\rho=\rho_A\otimes\rho_{BC}$ with $\rho_A$ a density matrix
on $\cH_A$ and $\rho_{BC}$ a density matrix on $BC$, then the partial transpose is given by
\[
\rho^\Gamma = \rho_A^T\otimes\rho_{BC},
\]
which is again a density matrix because the transpose of a positive semidefinite matrix is again positive semidefinite.

For general $\rho$, $\rho^\Gamma$ need no longer be positive definite; such quantum states are called
NPT states (for negative partial transpose). The class of NPT states is an important subset of the set of so-called
entangled states, which are the quantum states that have the greatest applicability in quantum information.
It is important to be able to quantify this entanglement. A very simple measure for doing so (although it does not reveal
all possible entangled states) is the so-called \textit{negativity of entanglement}, or simply \textit{negativity} of a quantum state,
first proposed in \cite{VW}.
Let us define $N(X):=||X||_1-\trace X$, the negativity function of a Hermitian matrix $X$.
Then the negativity of entanglement of $\rho$ between subsystems $A$ and $BC$ is defined as
\be
N_{A|BC}(\rho) = N(\rho^\Gamma).
\ee
Clearly, if $\rho^\Gamma$ is positive semidefinite, i.e.\ if $\rho$ is a state with positive partial transpose (PPT), then
its trace norm is equal to its trace, whence its negativity is zero.

One can also define the negativity of entanglement between subsystems $A$ and $B$ only, as
\be
N_{A|B}(\rho) = N(\trace_C \rho^\Gamma)
\ee
and between subsystems $A$ and $C$:
\be
N_{A|C}(\rho) = N(\trace_B \rho^\Gamma).
\ee
There are states for which $N_{A|B}$ is large and $N_{A|C}$ is small or even zero,
and there are states for which the opposite holds.
One can easily show (and we will do so below for pure states) that
both $N_{A|B}$ and $N_{A|C}$ are bounded above by $N_{A|BC}$; that is, the negativity of entanglement
can not increase under taking partial traces. These bounds are tight, as there exist states for which
$N_{A|B}=N_{A|BC}$ and there exist states for which $N_{A|C}=N_{A|BC}$.

We are now in the position to formulate the problem we wish to study in this paper, namely to prove the
widely held belief that $N_{A|B}$ and $N_{A|C}$ can not both be equal to $N_{A|BC}$.
This is the so-called \emph{monogamy of entanglement} property.
It has been proven in specific cases. The monogamy concept was introduced by Coffman, Kundu and Wootters \cite{CKW}, who proved it
for the smallest possible three-partite system, namely $d_A=d_B=d_C=2$ (three qubits), and for a different
measure of entanglement (the so-called concurrence). This result was then generalised to $n$-partite systems, each
subsystem still being 2-dimensional ($n$ qubits) by Osborne and Verstraete \cite{OV}.

Recently, He and Vidal \cite{HV}
conjectured that monogamy also holds in general 3-partite systems. Using negativity as an
entanglement measure they conjectured the following:
\begin{conjecture}[He-Vidal]
For any normalised complex vector $|\psi\rangle$ in the tensor product Hilbert space
$\cH=\cH_A\otimes \cH_B\otimes \cH_C$,
the following inequality holds:
\be
N^2(\trace_B |\psi\rangle\langle\psi|^\Gamma)
+ N^2(\trace_C |\psi\rangle\langle\psi|^\Gamma)
\le N^2(|\psi\rangle\langle\psi|^\Gamma).
\ee
\end{conjecture}
In the following section, we rephrase this problem as an inequality for certain block matrices that
can be readily understood without requiring any background knowledge in quantum information.
In the remainder of the paper we then
prove the conjecture in an important special case, using some well-established techniques of matrix analysis.
\section{Conversion to a block matrix problem\label{sec:convert}}
Given any set of orthonormal bases $\{|i\rangle\}$, $\{|j\rangle\}$ and $\{|k\rangle\}$
for the spaces $\cH_A$, $\cH_B$ and $\cH_C$, respectively,
we can write the pure state $|\psi\rangle$ we are considering as
\[
|\psi\rangle = \sum_{i=1}^{d_A} \sum_{j=1}^{d_B} \sum_{k=1}^{d_C} c_{ijk}|ijk\rangle.
\]
The coefficients $c_{ijk}$ can be rearranged into $d_A$ matrices $A_i$ with elements
\[
(A_i)_{jk} = c_{ijk}.
\]
We write $|A_i\rangle$ for the reshape of $A_i$ as a vector:
\[
|A_i\rangle = \sum_{jk} c_{ijk}|jk\rangle.
\]
Then $|\psi\rangle$ can be written in terms of the $|A_i\rangle$ as
\[
|\psi\rangle
= \sum_{ijk} c_{ijk} |i\rangle|jk\rangle
= \sum_i |i\rangle \otimes \sum_{jk} c_{ijk}|jk\rangle
= \sum_i |i\rangle \otimes |A_i\rangle,
\]
and
\[
\langle \psi| = \sum_i \langle i|\otimes \langle A_i|.
\]
The normalisation of $|\psi\rangle$ yields a condition on the $A_i$:
\be
1=\langle\psi|\psi\rangle = \sum_{i} \langle A_i|A_i\rangle = \sum_i\trace A_i^* A_i.\label{eq:normA}
\ee

\bigskip

The negativities can now be rewritten in terms of these matrices $A_i$.
For $N_{A|BC}$ we need $\rho^\Gamma$:
\beas
\rho^\Gamma&=&(|\psi\rangle\langle\psi|)^\Gamma \\
&=& \sum_{i,i'}(|i\rangle\langle i'|\otimes|A_i\rangle\langle A_{i'}|)^\Gamma \\
&=& \sum_{i,i'}|i'\rangle\langle i|\otimes|A_i\rangle\langle A_{i'}|.
\eeas
To find the modulus of this partial transpose, we first calculate the square.
\beas
\left(\rho^\Gamma\right)^2 &=&
\sum_{i,i',l,l'} \left(|i'\rangle\langle i|\otimes|A_i\rangle\langle A_{i'}|\right)\;
\left(|l'\rangle\langle l|\otimes|A_l\rangle\langle A_{l'}|\right) \\
&=& \sum_{i',l} |i'\rangle\langle l|\otimes \left(\sum_i |A_i\rangle\langle A_i|\right)\;\langle A_{i'}|A_l\rangle \\
&=& \left(\sum_{i',l} \langle A_{i'}|A_l\rangle\;\;|i'\rangle\langle l|\right)
\otimes \left(\sum_i |A_i\rangle\langle A_i|\right) \\
&=& \left(\left(\sum_{i'}|i'\rangle\langle A_{i'}|\right)\;\left(\sum_{l}|A_l\rangle\langle l|\right)\right)\;\;
\left(\left(\sum_{i}|A_i\rangle\langle i|\right)\;\left(\sum_{j}|j\rangle\langle A_{j}|\right)\right).
\eeas
Introducing the matrix
\be
\cA := \sum_i|A_i\rangle\langle i|, \label{eq:defcA}
\ee
which is a reshape of the vector of coefficients $c_{ijk}$, we can write $\left(\rho^\Gamma\right)^2$ as
\[
\left(\rho^\Gamma\right)^2 = \cA^*\cA \otimes \cA\cA^*.
\]
The modulus of $\rho^\Gamma$ is now simply the square root of this:
\[
|\rho^\Gamma| = |\cA|\otimes |\cA^*|.
\]
For the negativity we then get:
\beas
N_{A|BC} &=& ||\rho^\Gamma||_1-1 \\
&=& ||\cA||_1\;\;||\cA^*||_1 - 1
= ||\cA||_1^2 - 1 = ||\cA^*\cA||_{1/2}-1.
\eeas
It is worthwhile to point out that the matrix $\cA^*\cA$ is a Gram matrix,
as its elements in the $\{|i\rangle\}$ basis are given by
\[
\langle i|\cA^*\cA|i'\rangle = \langle A_i|A_{i'}\rangle.
\]
By the normalisation condition (\ref{eq:normA}), the trace of this Gram matrix is equal to $1$.

\bigskip

To find the other two negativities we need the partial traces of $\rho^\Gamma$:
\beas
\trace_C\rho^\Gamma
&=& \sum_{i,i'} |i'\rangle\langle i| \otimes \trace_C|A_i\rangle \langle A_{i'}| \\
&=& \sum_{i,i'} |i'\rangle\langle i| \otimes A_i A_{i'}^*
\eeas
and
\beas
\trace_B\rho^\Gamma
&=& \sum_{i,i'} |i'\rangle\langle i| \otimes \trace_N|A_i\rangle \langle A_{i'}| \\
&=& \sum_{i,i'} |i'\rangle\langle i| \otimes A_i^T (A_{i'}^*)^T \\
&=& \sum_{i,i'} |i'\rangle\langle i| \otimes \overline{A_i^* A_{i'}}.
\eeas
These partial traces can be expressed as $d_A\times d_A$ block matrices:
\[
\trace_C\rho^\Gamma = Z_1,\mbox{ and }
\trace_B\rho^\Gamma = \overline{Z_2},
\]
where
\be
Z_1:=
\left(
\begin{array}{ccc}
A_1A_1^* & A_2A_1^* & \ldots \\
A_1A_2^* & A_2A_2^* & \ldots \\
\vdots & \vdots &
\end{array}
\right) \label{eq:defZ1}
\ee
and
\be
Z_2:=
\left(
\begin{array}{ccc}
A_1^*A_1 & A_2^*A_1 & \ldots \\
A_1^*A_2 & A_2^*A_2 & \ldots \\
\vdots & \vdots &
\end{array}
\right). \label{eq:defZ2}
\ee
By the normalisation condition (\ref{eq:normA}), we have $\trace Z_1=\trace Z_2=1$.
The corresponding negativities are
\[
N_{A|B} = ||Z_1||_1-1\mbox{ and }
N_{A|C} = ||Z_2||_1-1.
\]

We can now reformulate the conjecture in terms of block matrices:
any $d_A$ block matrices $A_i$ satisfying the normalisation condition $\sum_i\trace A_i^*A_i=1$ also satisfy
the inequality
\be
\left(||Z_1||_1-1\right)^2 + \left(||Z_2||_1-1\right)^2 \le \left(||\cA^* \cA||_{1/2}-1\right)^2,
\label{ineq2}
\ee
where $\cA$, $Z_1$ and $Z_2$ are given by (\ref{eq:defcA}), (\ref{eq:defZ1}) and (\ref{eq:defZ2}), respectively.

\bigskip

It is possible to simplify the right-hand side of (\ref{ineq2}) by choosing a particular orthonormal basis
$\{|i\rangle\}$ for $\cH_A$.
Let the singular value decomposition of $\cA^*$ be given as $\cA^*=U\Sigma V^*$, where $U$ and $V$ are unitary matrices
of dimension $d_A$ and $d_Bd_C$, respectively, and $\Sigma$ is essentially diagonal.
If we choose $|i\rangle$ to be the $i$-th column of $U$, for all $i$, then
$\langle A_i|$ is the $i$-th row of $\Sigma V^*$. By this choice the vectors $|A_i\rangle$ become mutually orthogonal,
and the Gram matrix $\cA^* \cA$ becomes diagonal.
The right-hand side of (\ref{ineq2}) then simplifies by the identity
\[
|| \cA^* \cA||_{1/2}
= \left(\sum_i\sqrt{\langle A_i|A_i\rangle}\right)^2
= \left(\sum_i ||A_i||_2\right)^2.
\]

Inequality (\ref{ineq2}) is therefore equivalent to the somewhat simpler inequality
\be
\left(||Z_1||_1-1\right)^2 + \left(||Z_2||_1-1\right)^2 \le \left(\left(\sum_i ||A_i||_2\right)^2-1\right)^2,
\label{ineq3}
\ee
provided our choice of basis $\{|i\rangle\}$ is such that $\langle A_i|A_j\rangle=0$ for all $i\neq j$.

In pursuing a proof of (\ref{ineq2}) we may of course drop this condition and try and prove
(\ref{ineq3}) unconditionally, hoping that it is true in general.
Remarkably, however, inequalities (\ref{ineq3}) and (\ref{ineq2}) are equivalent even
without the orthogonality condition. This can be seen from the fact that replacing the positive semidefinite
matrix $\cA^*\cA$ by its diagonal can not decrease the $||\cdot||_{1/2}$ quasinorm.
In what follows we will focus on proving (\ref{ineq3}), unconditionally. 
By a further rescaling we can drop the normalisation condition $\sum_i\trace A_i^*A_i=1$, upon which (\ref{ineq3}) turns into
\bea
\left(||Z_1||_1-\trace Z_1\right)^2 + \left(||Z_2||_1-\trace Z_2\right)^2
&\le& \left(\left(\sum_i ||A_i||_2\right)^2-\sum_i ||A_i||_2^2\right)^2 \nonumber\\
&=& \left(\sum_{i\neq j} ||A_i||_2 ||A_j||_2\right)^2,
\label{ineq4}
\eea
which is the final form, as already advertised in the Introduction.
\section{Proof of a special case\label{sec:proof}}
The task of proving inequality (\ref{ineq4}) is a hard one because of the inequality's tightness.
It is easy to see that every term of the left-hand side of (\ref{ineq4}) is itself
bounded above by the right-hand side. In entanglement theory, this corresponds to the fact that the
negativity is a so-called entanglement monotone, which among other things means that it can not increase under
taking partial traces \cite{VW}. A matrix
analytical proof proceeds by first exploiting the triangle inequality to show that
$||Z_1||_1 \le \sum_{i,j} ||A_j A_i^*||_1$, and then
the Cauchy-Schwartz inequality to bound $\sum_{i,j} ||A_j A_i^*||_1$ by
$\sum_{i,j} ||A_j||_2 ||A_i^*||_2 = (\sum_i ||A_i||_2)^2$.

To prove (\ref{ineq4}), however, we must show that the sum of $\left(||Z_1||_1-\trace Z_1\right)^2$ and $\left(||Z_2||_1-\trace Z_2\right)^2$
is bounded above by the exact same expression that bounds each of the terms separately.
Finding the proof of that statement is an extremely delicate process, where picking up proportionality
constants has to be avoided at all costs.
Any such constant larger than 1 (no matter how close to 1) would ruin the tightness and render the result irrelevant.
For example, it is clear from the
above that (\ref{ineq4}) certainly holds with an extra factor of $2$ in the right-hand side (just add the inequalities
for each term separately) but this is a trivial result and says absolutely nothing about monogamy of negativity.

In what follows  we will restrict
to the case $d_A=2$ (i.e.\ system $A$ is a qubit); even this simple case already turned out to be a major undertaking.
To simplify notations, we will replace $A_1$ and $A_2$ by $A$ and $B$.
Then $Z_1$ and $Z_2$ are given by the $2\times 2$ block matrices
\[
Z_1 = \twomat{AA^*}{BA^*}{AB^*}{BB^*}\mbox{ and }
Z_2 = \twomat{A^*A}{B^*A}{A^*B}{B^*B}.
\]
Furthermore, we were obliged to restrict to the case $A=\id$. This requires taking $d_B=d_C$.
We will henceforth write $d$ for $d_B=d_C$.

In this case both terms of the left-hand side of (\ref{ineq4}) turn out to be less than one half the right-hand side.
Adding up then proves (\ref{ineq4}).
The goal is therefore to show, for all $d\times d$ matrices $B$,
\[
||Z_1||_1-\trace Z_1 \le \frac{1}{\sqrt{2}} \;2 \;||\id||_2 ||B||_2 = 2\sqrt{d/2} \; ||B||_2.
\]
Replacing $B$ by $B^*$ yields the corresponding inequality for $Z_2$.
Henceforth, we will write $Z$ for $Z_1$, and we have
\[
Z=\twomat{\id}{B}{B^*}{BB^*}.
\]
Noting that $||X||_1-\trace X=2\trace X_-$ for any Hermitian $X$, we can rewrite the inequality as
\be
\trace Z_- \le \sqrt{d/2} \; ||B||_2.
\label{ineqid}
\ee

Our proof proceeds by splitting this inequality into two inequalities.
First we show
\be
\trace Z_- \le \trace\sqrt{(BB^*-B^*B)_-}
\label{ineqid1}
\ee
and then we show
\be
\trace\sqrt{(BB^*-B^*B)_-}\le \sqrt{d/2} \; ||B||_2.
\label{ineqid2}
\ee
\subsection{Proof of inequality (\ref{ineqid1})}
It is well-known that any given square matrix $B$ is weakly unitarily equivalent to its Hermitian conjugate $B^*$.
Indeed, let $B=U|B|$ be the polar decomposition of $B$, then $B^*=|B|U^*$, so that
$B^*=U^* B U^*$. So, by multiplying $B$ on the left and on the right by some unitary matrix, we obtain $B^*$.
However, there is another way to relate $B$ and $B^*$ requiring only a left multiplication,
by extending both matrices.

Let $\Delta=BB^*-B^*B$ and let its Jordan decomposition be $\Delta=\Delta_+ - \Delta_-$, where $\Delta_\pm\ge0$.
Then
\[
BB^* + \Delta_- = B^*B+\Delta_+.
\]
By positive semidefiniteness of all four terms we can write this as
\[
\twovect{B}{\sqrt{\Delta_-}}\twovec{B^*}{\sqrt{\Delta_-}} =
\twovect{B^*}{\sqrt{\Delta_+}}\twovec{B}{\sqrt{\Delta_+}}.
\]
This immediately implies that there must exist a unitary matrix $U$ such that
\be
\twovec{B^*}{\sqrt{\Delta_-}} = U\twovec{B}{\sqrt{\Delta_+}}.\label{eq:BstarB}
\ee
These two block matrices are the abovementioned extensions of $B$ and $B^*$, respectively.
If $B$ is not square, it can be made so by zero-padding and the same statement therefore holds for general matrices $B$.

According to Cauchy's interlacing theorem,
the eigenvalues of an $m\times m$ principal submatrix $A'$ of an $n\times n$ Hermitian matrix $A$ satisfy the relation
$\lambda^\uparrow_j(A) \le \lambda^\uparrow_j(A')$ for $j=1,\ldots,m$
(there is also an upper bound, but we will not need it).
In particular, as $Z$ is a submatrix of the matrix
\[
Z_1 := \left(
\begin{array}{ccc}
\id & 0 & B \\
0 & \id & \sqrt{\Delta_+} \\
B^* & \sqrt{\Delta_+} & BB^*
\end{array}
\right)
\]
we have $\lambda_j^\uparrow (Z) \ge \lambda_j^\uparrow(Z_1)$ for $j=1,\ldots,2d$.

By (\ref{eq:BstarB}), and the fact
that for unitary $U$ two block matrices of the form
\[
\twomat{\id}{UX}{X^*U^*}{Y} \mbox{ and }
\twomat{\id}{X}{X^*}{Y},
\]
are equal up to a unitary conjugation and therefore
have the same spectrum, $Z_1$ has the same spectrum as
\[
Z_2 := \left(
\begin{array}{ccc}
\id & 0 & B^* \\
0 & \id & \sqrt{\Delta_-} \\
B & \sqrt{\Delta_-} & BB^*
\end{array}
\right).
\]
Now, $Z_2$ can be split as a sum of two matrices, the first one being positive semidefinite:
\[
Z_2 = Z_3+Z_4,\quad
Z_3:=\left(
\begin{array}{ccc}
\id & 0 & B^* \\
0 & \id & 0\\
B & 0 & BB^*
\end{array}
\right),\quad
Z_4:=
\left(
\begin{array}{ccc}
0 & 0 & 0 \\
0 & 0 & \sqrt{\Delta_-} \\
0 & \sqrt{\Delta_-} & 0
\end{array}
\right).
\]
By Weyl's monotonicity theorem,
we therefore have
\[
\lambda_j^\uparrow(Z_2) \ge \lambda_j^\uparrow(Z_4).
\]

The $d$ smallest eigenvalues of $Z_4$ are non-positive and given by $-\sqrt{\mu_j}$, where $\mu_j$ are the eigenvalues of $\Delta_-$.
Thus,
\[
\lambda_j^\uparrow (Z) \ge -\sqrt{\mu_j^\downarrow}, \mbox{ for }j=1,\ldots,d
\]
Furthermore, $Z_4$ has at most $d$ negative eigenvalues.
Tracing back through the previous argument then reveals that this is also true for $Z_2$, $Z_1$ and finally $Z$ itself.

So we have that the number of negative eigenvalues $n_-$ of $Z$ is at most $d$,
and they are larger than $-\sqrt{\mu_j}$.
Hence,
\[
\trace Z_- = \sum_{j=1}^{n_-} (-\lambda_j^\uparrow(Z)) \le \sum_{j=1}^{n_-} \sqrt{\mu_j^\downarrow} \le \sum_{j=1}^{d} \sqrt{\mu_j} = \trace\sqrt{\Delta_-},
\]
which is inequality (\ref{ineqid1}).
\subsection{Proof of inequality (\ref{ineqid2})\label{sec:ineqid2}}
In this section we prove that the inequality (\ref{ineqid2})
is valid for any $d\times d$ matrix $B$.
For convenience we will actually prove the equivalent statement that
\[
\trace\sqrt{(BB^*-B^*B)_+}\le \sqrt{d/2} \; ||B||_2;
\]
the latter turns into the former by replacing $B$ with $B^*$.

Note that $BB^*$ and $B^*B$ have the same eigenvalues, hence they are unitarily equivalent.
Another way to phrase the inequality is that
\[
\trace\sqrt{(L-ULU^*)_+}\le \sqrt{d/2} \; ||L||_2,
\]
for any unitary matrix $U$ and any non-negative diagonal matrix $L$.
One way to attack this problem is to first try and prove it for $U$ that are permutation matrices,
so that both $L$ and $ULU^*$ are diagonal, and then extend this result from the commutative case to the general case.
It turns out that this extension can indeed be done thanks to a theorem by Drury.

In \cite{drury2}, Drury stated the following theorem (without explicit proof, but with the remark that it can 
be proven easily using the method he has developed in a preceding publication, \cite{drury1}):
\begin{theorem}[Drury]
Let $X$ and $Y$ be $d\times d$ Hermitian matrices with given eigenvalues
$x_1\ge x_2\ge\cdots\ge x_d$ and $y_1\ge y_2\ge\cdots\ge y_d$, respectively. Let $I=[x_n+y_n,x_1+y_1]$.
Let the function $\phi:I\to\bR$ be isoclinally metaconvex on $I$. Then
\[
\trace \phi(X+Y) \le \max_{\pi\in S_d} \sum_{j=1}^d \phi(x_j+y_{\pi(j)}).
\]
\end{theorem}
The class of isoclinally metaconvex (IM) functions has been introduced by Drury in \cite{drury2} exactly for this purpose:
\begin{definition}
Let $I$ be an interval in $\bR$. An infinitely differentiable function $\phi:I\to\bR$ is said to be
IM on $I$ if whenever $t_1,t_2\in I$ with $t_1\neq t_2$ and
$\phi'(t_1)=\phi'(t_2)$, then $\phi''(t_1)+\phi''(t_2)>0$.
\end{definition}
For example, strictly concave and strictly convex functions are both IM.
It is possible for other functions to be in this class as well, provided that for every point where the curvature
is negative there is another point with the same gradient and with positive curvature greater in absolute value.

This theorem would allow us to reduce the problem of proving inequality (\ref{ineqid2}) to the commutative case
if only the function $x\mapsto f(x)=\sqrt{x_+}$ were IM. Clearly it is not, as it is
not even differentiable. However, $f(x)$ can be approximated arbitrarily well by a sequence of
IM functions, as shown in the appendix, and this is all what is needed.

Hence, Drury's result when applied to the matrices $X=BB^*$ and $Y=-B^*B$ implies that inequality (\ref{ineqid2}) is valid
if we can show that the inequality
\[
\left(\sum_{i=1}^d \sqrt{(\mu_i-\mu_{\pi(i)})_+}\right)^2 \le \frac{d}{2}\;\sum_{i=1}^d \mu_i
\]
holds for any permutation $\pi\in S_d$, and for any set $\mu_i$ of non-negative numbers (the eigenvalues of $BB^*$).
Without loss of generality we can assume that $\mu_1\ge\mu_2\ge\cdots\ge\mu_d$ and $\sum_i\mu_i=1$.

%
The key to the proof is to decompose a given permutation $\pi\in S_d$ in what we will call here \emph{maximal ascending chains}
(MA chains).
Let an ascending chain be a sequence of increasing integers from $\{1,2,\ldots,d\}$ such that
the image under $\pi$ of each integer
in the chain is given by the next integer in the chain (if any). That is, it is a sequence $I:=(i_1,i_2,\ldots,i_r)$
such that $i_{j+1}>i_j$ and $i_{j+1}=\pi(i_j)$, for $j=1,2,\ldots,r-1$.
An MA chain is one that is as long as possible. For a general permutation, more than one such chain
may exist. Clearly, chains are disjoint.

For example, the permutation $\twovec{1234}{2341}$ has one MA chain, namely
$I=(1,2,3)$. The element 4 is not included because its image is 1, which is less than 4.
The permutation $\twovec{1234}{3421}$ has two such chains, namely $I_1=(1,3)$ and $I_2=(2,4)$.

To proceed with the proof, we split the sum $\sum_i\sqrt{(\mu_i-\mu_{\pi(i)})_+}$ into several components,
one per MA chain of the permutation $\pi$.
Let the lengths of the various MA chains $I_1$, $I_2$, \ldots, $I_K$ of a permutation
be $r_1,r_2,\ldots,r_K$, respectively.
Clearly, as MA chains are disjoint, the $r_k$ sum up to at most $d$.
Then we split the sum as follows:
\[
\sum_{i=1}^d \sqrt{(\mu_i-\mu_{\pi(i)})_+}
= \sum_{k=1}^K \sum_{i\in I_k} \sqrt{(\mu_i-\mu_{\pi(i)})_+}.
\]
We can do this because the $i$-th term has a nonzero contribution to the sum unless $i$ appears in some MA chain.
Indeed, if $i$ does not appear in any of the MA chains, this means that $\pi(i)<i$,
whence, by the ordering of the $\mu_i$, we have $\mu_{\pi(i)}>\mu_i$, so that $(\mu_i-\mu_{\pi(i)})_+=0$.

Let us now consider one such component, for a chain $I=(i_1,i_2,\ldots,i_r)$ of length $r$:
\[
\sum_{j=1}^{r-1} \sqrt{(\mu_{i_j}-\mu_{\pi(i_j)})_+}.
\]
Because the $i_j$ form an MA chain, we have $\mu_{i_j}>\mu_{\pi(i_j)}$ and $\pi(i_j)=i_{j+1}$.
We can therefore simplify this sum as
\[
\sum_{j=1}^{r-1} \sqrt{\mu_{i_j}-\mu_{i_{j+1}}}
\]
(for example, the chain $(1,2,3)$ mentioned before corresponds to the component
$\sqrt{\mu_1-\mu_2}+\sqrt{\mu_2-\mu_3}$).
We now claim that this sum is bounded above by the quantity
\[
\left(\frac{r}{2}\sum_{j=1}^r \mu_{i_j}\right)^{1/2}
\]
(in the example, by $\sqrt{(3/2)(\mu_1+\mu_2+\mu_3)}$).

For $r=2$ this is trivially true, as the sum has only one term:
\[
\left(\sum_{j=1}^{r-1} \sqrt{\mu_{i_j}-\mu_{i_{j+1}}}\right)^2
= \mu_{i_1}-\mu_{i_2} \le \mu_{i_1}+\mu_{i_2}.
\]
For $r>2$ we can exploit the following inequality, which can be seen as a H\"older-type inequality
for the $l_{1/2}$-(quasi)-norm: for any vector $x$ with non-negative real elements $x_j$,
and any probability vector $p$ (that is, $p_j\ge 0$ and $\sum_j p_j=1$),
\be
\left(\sum_{j=1}^d \sqrt{x_j}\right)^2 \le \sum_{j=1}^d \frac{x_j}{p_j}.
\ee
We will apply this in the following instance: $d=r-1$, $x_j=\mu_{i_j}-\mu_{i_{j+1}}$ and
$p_1=2/r$ and $p_2=\cdots=p_{r-1}=1/r$, to obtain, as required,
\beas
\left(\sum_{j=1}^{r-1} \sqrt{\mu_{i_j}-\mu_{i_{j+1}}}\right)^2
&\le& \frac{r}{2}(\mu_{i_1}-\mu_{i_2}) + r(\mu_{i_2}-\mu_{i_3})+\cdots+r(\mu_{i_{r-1}}-\mu_{i_r}) \\
&=& \frac{r}{2}(\mu_{i_1}+\mu_{i_2})-r\mu_{i_r}
\le \frac{r}{2}\sum_{j=1}^r \mu_{i_j}.
\eeas

Having  one such bound per MA-chain component, we can now easily
get a bound on the entire sum. By the previous result we have
\[
\sum_{i=1}^d \sqrt{(\mu_i-\mu_{\pi(i)})_+}
\le \sum_{k=1}^K \sqrt{\frac{r_k}{2}\sum_{i\in I_k} \mu_i}.
\]
We can now simply exploit the Cauchy-Schwarz inequality and find the upper bound
\[
\sum_{k=1}^K \sqrt{\frac{r_k}{2}\sum_{i\in I_k} \mu_i}
\le \left(\sum_{k=1}^K \frac{r_k}{2}\right)^{1/2} \; \left(\sum_{k=1}^K \sum_{i\in I_k} \mu_i \right)^{1/2}
\le \sqrt{\frac{d}{2}}\;\sqrt{\sum_{i=1}^d\mu_i},
\]
which ends the proof.
\section{Appendix}
Here we prove the statement used in Section \ref{sec:ineqid2} that the function $f(x)=\sqrt{x_+}$ can be approximated
arbitrarily well by IM functions. More precisely, we show that there exists a sequence of IM functions that
converges uniformly to $f(x)$.
Many functions do so, but we construct this sequence in such a way that its
metaconvexity is easy to prove.

We start by defining a particular function $h(x)$ and then show two things: first, that $h(x)$ is IM and second,
that $|h(x)-f(x)|$ is bounded by a finite constant $c>0$.
Using such a function $h(x)$ we can easily construct a sequence of IM functions converging uniformly to $f(x)$:
we just have to consider the functions $h_s(x):=h(sx)/\sqrt{s}$. These functions inherit the property of being IM from $h(x)$,
and $|h_s(x)-f(x)| = |h(sx)-f(sx)|/\sqrt{s}<c/\sqrt{s}$, which tends to 0 as $s$ tends to $+\infty$, proving their
uniform convergence.

To construct $h(x)$ consider the functions
$w(x)=(1/2)(x^2+1)^{-1/4}$ and $\alpha(x) = 1+\exp(-x)$,
and let $g(x)=w(\alpha(x)x)$. The function $\alpha(x)$ satisfies $\alpha(x)\ge 1$, is monotonically decreasing, and
tends to 1 as $x$ tends to $+\infty$.
Our function $h(x)$ of choice is the integral of $g(x)$, namely $h(x)=\int_{-\infty}^x dy\;g(y)$.
Note that for $x$ tending to $+\infty$, $w(x)$ tends to $1/(2\sqrt{x})$, so that in that regime $h(x)$ tends
to $\sqrt{x}$ plus some finite constant arising from the integration over all smaller values of $x$.

We first show that $h(x)$ is IM.
This involves the first and second derivatives of $h$, which are given by:
\beas
h'(x) = g(x) &=& w(\alpha(x)x) \\
h''(x) = g'(x) &=& w'(\alpha(x)x) \; (\alpha'(x)x+\alpha(x)).
\eeas
We therefore need to show that distinct $x_1$ and $x_2$ with the same value of $g(x)$ must satisfy $g'(x_1)+g'(x_2)>0$.
It is essential that $w(x)$ is an even function that is monotonically increasing for $x<0$, and
monotonically decreasing for $x>0$, so that any pair of distinct $x$ having the same $g(x)$ must have opposite sign.
Let $x_1<0$ and $x_2>0$ be such points. By the evenness of $w(x)$,
this is so if and only if $-\alpha(x_1)x_1 = \alpha(x_2)x_2$.
For such points, the factor $w'(\alpha(x)x)$ in $g'(x)$ has the same absolute value (again by virtue of $w$ being even),
and is positive for $x_1$ and negative for $x_2$. The condition $g'(x_1)+g'(x_2)>0$ is therefore equivalent to
\[
(\alpha'(x_1)x_1+\alpha(x_1)) - (\alpha'(x_2)x_2+\alpha(x_2))>0.
\]
This condition is easily seen to be satisfied as $\alpha'(x)x+\alpha(x)=-\theta x \exp(-\theta x)+1+\exp(-\theta x)$
is always larger than 2 for $x<0$ and less than 2 for $x>0$.
This proves that $h(x)$ is IM.

Secondly, we have to show that $h(x)$ is an approximation of $f(x)=\sqrt{x_+}$, in the sense that
$|h(x)-f(x)|$ is bounded by a finite constant.

For $x<0$ we have $f(x)=0$ and $h(x)>0$. To show that $h(x)$ is bounded above for $x<0$ we only have to show that $h(0)$ is
finite, since $h(x)$ is an increasing function (as $w(x)>0$).
Since $w(x)<1/(2\sqrt{-x})$ for $x<0$ and $\alpha(x)>\exp(-\theta x)$, we get, indeed,
\[
h(0) = \int_{-\infty}^0 dy \; w(\alpha(y) y)
< \int_{-\infty}^0 dy \; \frac{1}{2\sqrt{\exp(-\theta y)y}}
= \sqrt{\frac{\pi}{2\theta}}.
\]

For $x>0$, $f(x)=\sqrt{x}$. As $\alpha(x)>1$, we have that $h'(x) = w(\alpha(x)x) < w(x)$. For $x>0$,
we also have $w(x)<1/(2\sqrt{x})=f'(x)$, so that $h'(x)<f'(x)$. Integrating over $x$ yields
$h(x)-h(0)\le f(x)-f(0)$ from which we obtain the upper bound $h(x)-f(x)<h(0)-f(0) = h(0)$, which is finite.

To obtain a lower bound we can exploit the two inequalities
\[
w(x)=\frac{1}{2(1+x^2)^{1/4}} > \frac{1}{2\sqrt{x}} - \frac{1}{8x^{5/2}}
\mbox{ and }
\frac{1}{\sqrt{1+\exp(-x)}} > 1-\frac{1}{2}\exp(-x).
\]
This yields
\beas
h'(x) = w(\alpha(x)x) &>& \frac{1}{2\sqrt{x}\sqrt{1+\exp(-x)}} - \frac{1}{8(1+\exp(-x))^{5/2}} \\
&>& \frac{1}{2\sqrt{x}}\left(1-\frac{1}{2}\exp(-x)\right) - \frac{1}{8x^{5/2}}
\eeas
so that
\[
h'(x)-f'(x)> -\frac{\exp(-x)}{4\sqrt{x}} - \frac{1}{8x^{5/2}}.
\]
Integrating from 1 to $x$ yields, for $x>1$,
\[
h(x)-f(x)>h(1)-\left(\frac{1}{4}\int_1^x dx \;\frac{\exp(-x)}{\sqrt{x}} + \frac{1}{8}\int_1^x dx\;\frac{1}{x^{5/2}}\right).
\]
The first integral is bounded above by $\int_0^\infty dx \;\exp(-x)/\sqrt{x} = \sqrt{\pi}$
and the second integral is equal to $(2/3)(1-x^{-3/2})$, which is bounded above by $2/3$.
Thus, for $x>1$, $h(x)-f(x)$ is bounded below by a finite constant.
It is clear that, for $0<x<1$, $h(x)-f(x)$ is bounded below as well since $h(x)>0$ and $f(x)<1$.
We conclude that $|h(x)-f(x)|$ is bounded everywhere by a finite constant.
\section*{Acknowledgments}
This work is supported by an Odysseus grant from the Flemish FWO.
Thanks to Frank Verstraete for discussions and rekindling my interest in entanglement theory.


\begin{thebibliography}{9}
\bibitem{cargese} K.M.R.~Audenaert, ``Mathematical aspects of quantum information theory'',
in: Physics and Theoretical Computer Science, J.-P.~Gazeau et al. (Eds.), 3--24, IOS Press (2007).
\bibitem{AK} K.M.R.~Audenaert and F.~Kittaneh, ``Problems and Conjectures in Matrix and Operator Inequalities'',
in: Operator Theory, J.~Zemanek ed., Banach Center Publications Series. (In Press). See also eprint arXiv:1201.5232.
\bibitem{CKW} V.~Coffman, J.~Kundu and W.~Wootters, ``Distributed entanglement'', Phys.\ Rev.\ A \textbf{61}, 052306 (2000).
\bibitem{drury1} S.W.~Drury, ``On symmetric functions of the eigenvalues of the sum of two Hermitian matrices'',
Linear Algebra Appl.\ \textbf{176}, 211-222 (1992).
\bibitem{drury2} S.W.~Drury, ``Maximizing traces of matrix functions'',
Linear Algebra Appl.\ \textbf{387}, 221--234 (2004).
\bibitem{HV} H.~He and G.~Vidal, ``Disentangling theorem and monogamy for entanglement negativity'',
eprint arXiv:1401.5843 (2014).
\bibitem{OV} T.J.~Osborne and F.~Verstraete, ``General monogamy inequality for bipartite qubit entanglement'',
Phys.\ Rev.\ Lett.\ \textbf{96}, 220503 (2006).
\bibitem{VW} G.~Vidal and R.F.~Werner, ``A computable measure of entanglement'', Phys.\ Rev. A \textbf{65}, 032314 (2002).
\end{thebibliography}
\end{document}